\documentclass{elsarticle}
\usepackage{epsfig}
\usepackage{float}
\usepackage{amssymb}
\usepackage{mathrsfs}
\usepackage{enumerate}
\usepackage{hhline}
\usepackage{amsmath}
\usepackage{indentfirst}

\def\H{\mathcal{H}}

\usepackage{graphicx}

\newtheorem{proposition}{Proposition}[section]

\newtheorem{lemma}[proposition]{Lemma}

\newtheorem{theorem}[proposition]{Theorem}

\newtheorem{corollary}[proposition]{Corollary}

\newtheorem{definition}[proposition]{Definition}

\def\p{\noindent{\bf Proof. }}
\def\q{\hspace*{\fill}$\Box$\medskip}
\usepackage{color}
\usepackage{cite}
\usepackage{subfigure}

\journal{ }
\begin{document}
	\begin{frontmatter}
		\title{The effect on the spectral radius of $r$-graphs by grafting or contracting edges}

		\author[t3]{Wei Li}
		\address[t3]{School of Computer and Information Science, Fujian Agriculture and Forestry University,  Fuzhou, Fujian, 350002, P. R. China}
		\address[t2]{Center for Discrete Mathematics and Theoretical Computer Science, Fuzhou University, Fuzhou, Fujian, 350003, P. R. China}
		
		\author[t2]{An Chang}
		
		\begin{abstract}
			Let $\mathcal{H}^{(r)}_n$ be the set of all connected $r$-graphs with given size $n$. In this paper, we investigate the effect on the spectral radius of $r$-uniform hypergraphs by grafting or contracting an edge and then give the ordering of the $r$-graphs with small spectral radius over $\H^{(r)}_n$, when $n\geq 20$.	
		\end{abstract}

		
	\end{frontmatter}

	\section{Introduction}
\begin{figure}[htbp]
		\centering
		\includegraphics[width=0.8\textwidth]{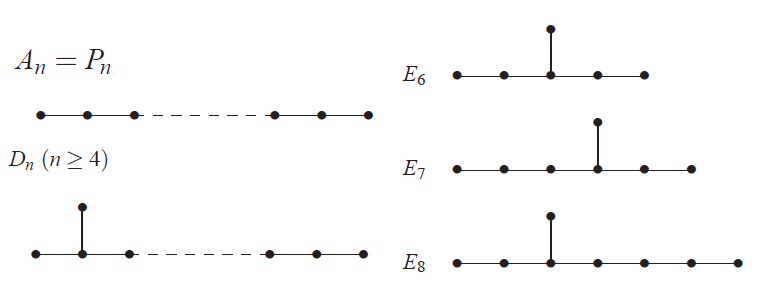}
		\caption{The graphs with spectral radius less than 2}
		\label{fig2}
	\end{figure}

	As we know, the spectral radius of a graph $G$, denoted by $\rho(G)$, is the largest eigenvalue of its adjacency matrix. It plays a very important role in the spectra graph theory. In 1970, Smith \citep{Smith1970Some} determined all connected graphs with spectral radius at most 2, in which those graphs with spectral radius less than 2 are called Dynkin Diagram (see Fig.\ref{fig2}).
	
	In 2005,  the eigenvalues of higher-order tensors were introduced by Qi\citep{qi2005} and Lim \citep{lim2005} independently. Since then, the study on the spectra of tensors and their various applications have been attracted much attention and interest. In addition, there are a lot of results concerning the spectral theory of uniform hypergraphs via tensors. In 2009, Bul$\grave{o}$ and Pelillo \citep{pelillo2009} gave new bounds on the clique number of a graph based on analysis of the largest eigenvalue of the adjacency tensor of a uniform hypergraph. In 2012, Cooper and Dulte \citep{cooper} analyzed eigenvalues of the adjacency tensor of a uniform hypergraph, and proved a number of natural analogs of basic resutls in spectral graph theory. In 2015, Li, Shao and Qi \citep{Lihh} studied some extremal spectral properties of the classes of uniform supertrees with $n$ vertices, and determined that the hyperstar attains uniquely the maximal spectral radius among all uniform supertrees with given size. Yuan et al. \citep{yuan2015, yuan2016} determined the ordering of supertrees with larger spectral radius. Fan et al.\citep{fanyz2015} determined the hypergraphs with maximum spectral radius over all unicyclic hypergraphs, over linear or power unicyclic hypergraphs with given girth, over linear or power bicyclic hypergrphs, respectively. Some other study in the spectra of uniform hypergraphs can be found in \citep{shao,xie2012,xie2013,xieLAA,pearson,kchang2008,kchang2011,friedland2013}.
	
	In 2014, Lu and Man \citep{Lu2014Connected} extended the Smith's result to $r$-uniform hypergraphs and described all connected $r$-uniform hypergraphs with spectral radius at most $(r-1)!\sqrt[r]{4}$. Futhermore, to approximate the spectral radius of uniform hypergraphs,  they explored a new technique using weighted incidence matrix. In \citep{lM2}, Lu and Man also determined the $r$-uniform hypergraph with spectral rasius at most $2\sqrt[r]{2+\sqrt{5}}$.

    Let $\mathcal{H}^{(r)}_n$ be the set of all connected $r$-graphs with given size $n$. In this paper, we investigate the effect on the spectral radius of $r$-uniform hypergraphs by grafting or contracting an edge. Based on these results, we compare the sepctral radius of $r$-uniform hypergraphs classified in \citep{Lu2014Connected} and gave the ordering of $r$-graphs with small spectral radius over $\H^{(r)}_n$, when $n\geq 20$.
    
    This paper is organized as follows. In section 2, the notation and some important lemmas are listed. In section 3, we study the effect on the spectral radius of $r$-graphs under perturbations and then the ordering is given in  section 4.
	
\section{Preliminary}
 A {\it hypergraph} $H$ is a pair $(V,E)$. The elements of $V=V(H)=\{v_1,v_2,\cdots,$  $v_{\nu}\}$ are referred to as {\it vertices} and the elements of $E=E(H)=\{e_1,e_2,\cdots,$ $e_n\}$ are called {\it edges}, where $e_i\subseteq V$ for $i\in[n]$. A hypergraph $H$ is said to be {\it r-uniform} for an integer $r\geq2$, if for all $e_i\in E(H)$, $|e_i|=r$, where $i\in[n]$. Throughout this paper, $n$ always denotes the number of edges in $H$ and we often use the term {\it $r$-graph} in  place of $r$-uniform hypergraph for short. Obviously, $2$-graph is the general graph we usually say. A hypergraph $H$ is called {\it simple} if every pair of edges intersects at most one vertex. In fact, any non-simple hypergraph contains a cycle: $v_1e_1v_2e_2v_1$, {\it i.e.}, $v_1,v_2\in e_1\cap e_2$. A vertex with degree one is called a {\it leaf} vertex.

Let $H$ be an $r$-graph with $\nu$ vertices. Then a vector $x=(x_1,x_2,\cdots,x_\nu)\in \mathbb{R}^\nu$  can be considered as a function from $V(H)$ to the real number set, where each vertex $v_i$ is mapped to $x_i$.  Naturally, we can say $x_i$ is the value of the vertex $v_i$, denoted by $x(v_i)$.

Considering the polynomial form from $x=(x_1,x_2,\cdots,x_\nu)\in \mathbb{R}^\nu$ to a real number, which is defined by

 $$P_H(x)=r!\sum_{\{v_{i_1},v_{i_2},\cdots v_{i_r}\}\in E(H)}x(v_{i_1})x(v_{i_2})\cdots x(v_{i_r}),$$
the spectral radius of $H$, denoted by $\rho(H)$, is the maximum value of $P_H(x)$ over the $r$-norm sphere. In \citep{kchang2008, friedland2013, yang2010}, the Perron-Frobenius Theorem deduces the following result concerning with the spectral radius and the corresponding eigenvector for a connected $r$-graph.

\begin{lemma}\citep{kchang2008, friedland2013, yang2010}
	If $H$ is a connected r-uniform hypergraph, then there exists a nonnegative eigenvector corresponding  to $\rho(H)$. Moreover, this nonnegative vector is unique and called by {\it Perron-Frobenius vector}.
\end{lemma}

In what follows, some useful definitions and results proposed by Lu and Man \citep{Lu2014Connected} are listed.
\begin{definition}\citep{Lu2014Connected}
	A {\it weighted incident matrix} $B$ of a hypergraph $H$ is a $|V|\times |E|$ matrix such that for any vertex $v$ and edge $e$, the entry $B(v,e)>0$ if $v\in e$ and $B(v,e)=0$ if $v\notin e$.
\end{definition}
\begin{definition}\citep{Lu2014Connected}\label{normal}
	A hypergraph $H$ is called {\it $\alpha$-normal} if there exists a weighted incidence matrix $B$ satisfying
	\begin{enumerate}
		\item  $\sum_{e: v\in e}B(v,e)=1$, for any $v\in V(H)$.
		\item  $\prod_{v\in e} B(v,e)=\alpha $, for any $e\in E(H)$.
	\end{enumerate}
	
	Moreover, the incidence matrix $B$ is called {\it consistent} if for any cycle $v_0e_1v_1e_2$ $\cdots e_lv_l$, where $(v_l=v_0)$, $$ \prod_{i=1}^{l}\frac{B(v_i,e_i)}{B(v_{i-1},e_{i})}=1.$$
	In this case, we call $H$ {\it consistently $\alpha$-normal}.
\end{definition}
	 
\begin{lemma}\citep{Lu2014Connected}\label{lem1}
	Let $H$ be a connected r-uniform hypergraph. 
$H$ is consistently $\alpha$-normal if and only if $\alpha=((r-1)!/\rho(H))^r$.
\end{lemma}
{\bf Remark.}  If $H$ is consistently $\alpha$-normal and $x$ is the Perron-Frobenius vector, then for any edge $e=\{v_{i_1},v_{i_2},\cdots, v_{i_r}\}\in E(H)$, 
\begin{equation}\label{eq1}
B(v_{i_1},e)^{1/r}x(v_{i_1})= B(v_{i_2},e)^{1/r}x(v_{i_2})=\cdots =B(v_{i_r},e)^{1/r}x(v_{i_r}).
\end{equation}
\begin{definition}\citep{Lu2014Connected}\label{super-d}
	A hypergraph  $H$ is called {\it $\alpha$-supernormal} if there exists a weighted incidence matrix $B$ satisfying
	\begin{enumerate}
		\item  $\sum_{e:\; v\in e}B(v,e)\geq 1$, for any $v\in V(H)$.
		\item  $\prod_{v\in e} B(v,e)\leq \alpha $, for any $e\in E(H)$.
	\end{enumerate}	
	Moreover, $H$ is called {\it strictly $\alpha$-supernormal} if it is $\alpha$-supernormal but not $\alpha$-normal.
\end{definition}	

\begin{lemma}\citep{Lu2014Connected}\label{super-t}
	Let $H$ be an r-uniform hypergraph. If $H$ is strictly and consistently $\alpha$-supernormal, then the spectral radius of $H$ satisfies
	$$\rho(H)>(r-1)! \alpha^{-\frac{1}{r}}.$$
	
\end{lemma}

\begin{definition}\citep{Lu2014Connected}\label{sub-d}
	A hypergraph  $H$ is called {\it $\alpha$-subnormal} if there exists a weighted incidence matrix $B$ satisfying
	\begin{enumerate}
		\item  $\sum_{e:\; v\in e}B(v,e)\leq 1$, for any $v\in V(H)$.
		\item  $\prod_{v\in e} B(v,e)\geq \alpha $, for any $e\in E(H)$.
	\end{enumerate}	
	Moreover, $H$ is called {\it strictly $\alpha$-subnormal} if it is $\alpha$-subnormal but not $\alpha$-normal.
\end{definition}	

\begin{lemma}\citep{Lu2014Connected}\label{sub-t}
	Let $H$ be an r-uniform hypergraph. If $H$ is $\alpha$-subnormal, then the spectral radius of $H$ satisfies
	$$\rho(H)\leq(r-1)! \alpha^{-\frac{1}{r}}.$$
	Moreover, if $H$ is strictly $\alpha$-subnormal then $\rho(H)<(r-1)!\alpha^{-\frac{1}{r}}$.
\end{lemma}

 $H'=(V',E')$ is said to be {\it extending} from $H=(V,E)$, if $H'=(V',E')$ is obtained by adding a new vertex in each edge of a $r$-uniform hypergraph $H$. Then $H'$ is $(r+1)$-uniform and  one can also say $H$ extends $H'$.  $A_n^{(3)}$ (see Fig. \ref{small}) is the 3-graph extending from $A_n$ (see Fig. \ref{fig2}).  $A_n^{(r)}$ is an $r$-graph extending from $A_n$ by $r-2$ times, called a simple path (or path for short) naturally.
 
 Observe that in any $\alpha$-normal incident matrix $B$, if an edge $e$ is incident with a leaf vertex $v$, then $B(v,e)=1$. We will omit this trivial value throughout this article. Moreover, it also leads to following lemma. 
 \begin{lemma}\label{same}\citep{Lu2014Connected}
 	If $H$ extends $H'$, then $H$ is consistently $\alpha$-normal if and only if $H'$ is consistently $\alpha$-normal for the same value of $\alpha$.
 \end{lemma}
	
\begin{figure}[htbp]
	\begin{minipage}[h]{0.5\textwidth}
	\centering
	\includegraphics[width=1.0\textwidth]{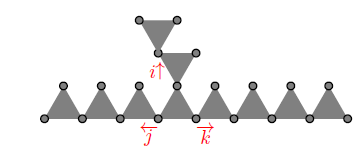}
	\caption{$F^{(3)}_{i,j,k}$}
	\label{fijk}
\end{minipage}
\begin{minipage}[htbp]{0.5\textwidth}
	\centering
	\includegraphics[width=1.0\textwidth]{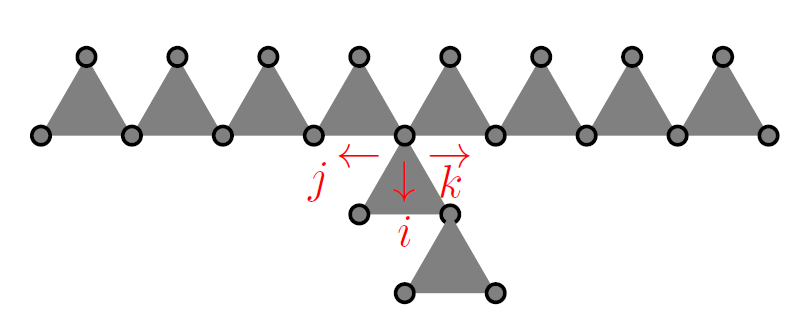}
	\caption{$E^{(3)}_{i,j,k}$}
	\label{eijk}
\end{minipage}
\end{figure}

\begin{figure}[htbp]
	\centering
	\includegraphics[width=0.8\textwidth]{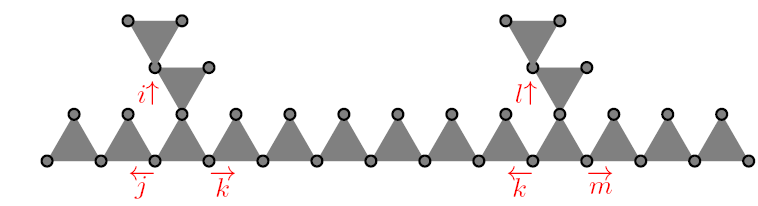}
	\caption{$G^{(3)}_{i,j:k:l,m}$}
	\label{gijlm}
\end{figure}

We denote by $E^{(3)}_{i,j,k}$ the 3-graphs obtained by attaching three paths of length $i,j,k$ to one vertex (see Fig.\ref{eijk}).
Denoted by $F_{i,j,k}^{(3)}$ the 3-graphs obtained by attaching three paths of length $i,j,k$ to each vertex of one edge (see Fig.\ref{fijk}). Denoted by $G^{(3)}_{i,j:k:l,m}$ the 3-graph obtained by attaching four paths of length $i,j,l,m$ to four ending vertices of path of length $k+2$ (see Fig.\ref{gijlm}).
\begin{figure}[htbp]
	\centering
	\includegraphics[width=1.0\textwidth]{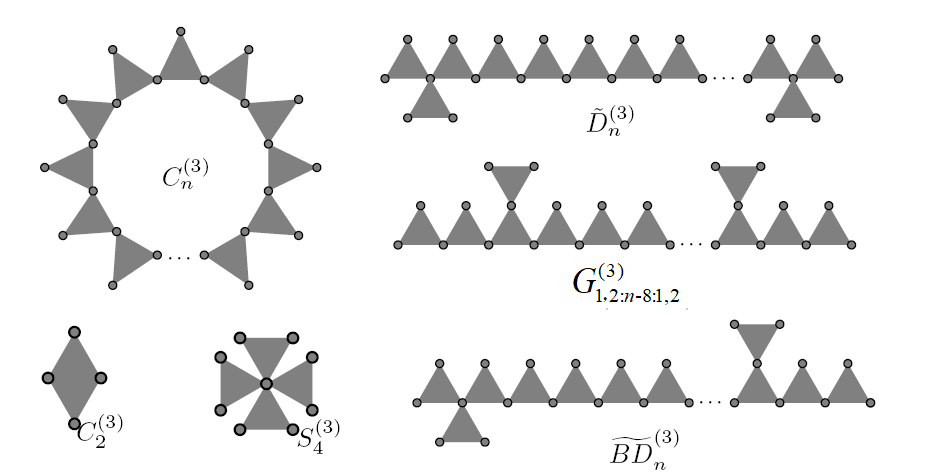}
	\caption{Some 3-graphs with spectral radius $2\sqrt[3]{4}$ in Theorem \ref{3-uniform-e}}
	\label{equal}
\end{figure}
\begin{figure}[htbp]
		\centering
		\includegraphics[width=1.0\textwidth]{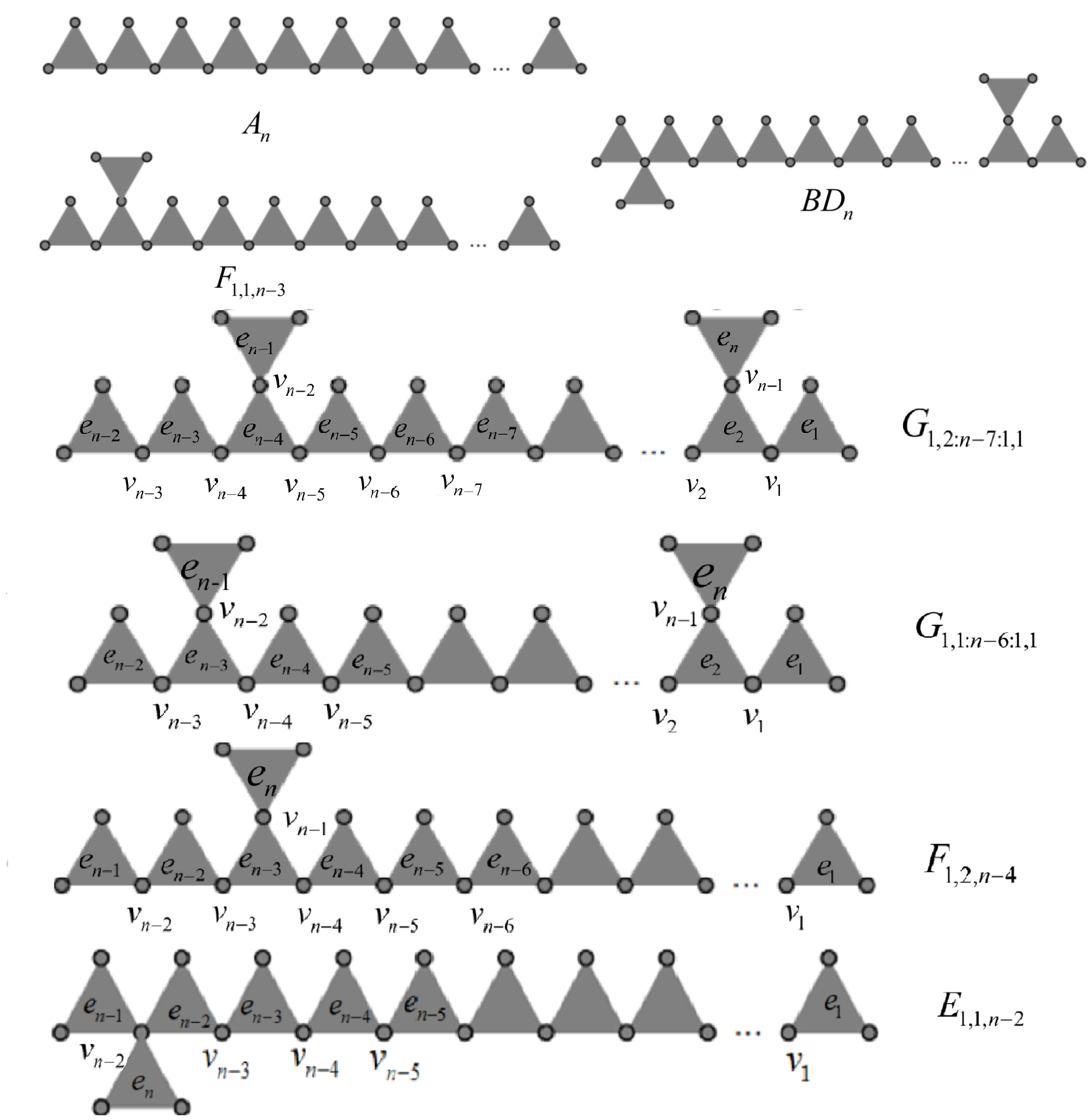}
		\caption{Some 3-graphs with small spectral radius in Theorem \ref{r-uniform-s}}
		\label{small}
	\end{figure}
\begin{figure}[htbp]
	\centering
	\includegraphics[width=0.9\textwidth]{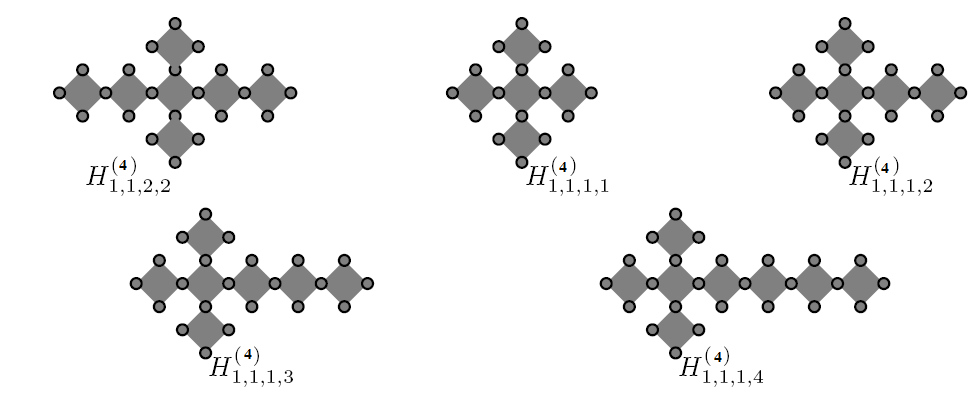}
	\caption{Some 4-graphs with spectral radius at most $3!\sqrt[4]{4}$}
	\label{4-uniform}
\end{figure}

\begin{theorem}\citep{Lu2014Connected}\label{r-uniform-s}
Let $r\geq 4$. If the spectral radius of a connected r-uniform hypergraph $H$ is less than $(r-1)!\sqrt[r]{4}$, then $H$ must be one of the following graphs:
\begin{enumerate}
	\item[(1)] $A_n^{(r)}$, $E_{1,1,n-2}^{(r)}$, $F_{1,1,n-3}^{(r)}$, $F_{1,2,n-4}^{(r)}$, $G_{1,1:n-6:1,1}^{(r)}$, $G_{1,2:n-7:1,1}^{(r)}$, $BD_n^{(r)}$, which are extending from the hypergraphs shown in Fig.6 by $r-3$ times.
	\item[(2)]  $E_{1,2,2}^{(r)}$, $E_{1,2,3}^{(r)}$, $E_{1,2,4}^{(r)}$, $F_{2,3,3}^{(r)}$, $F_{2,2,j}^{(r)}$ (for $2\leq j\leq 6$), $F_{1,3,j}^{(r)}$ (for $3\leq j\leq 13$), $F_{1,4,j}^{(r)}$ (for $4\leq j\leq 7$), $F_{1,5,5}^{(r)}$ and $G_{1,1:j:1,3}^{(r)}$ (for $0\leq j\leq 5$), which are the $r$-graphs extending from those corresponding 3-graphs by $r-3$ times.
	\item[(3)] $H_{1,1,1,1}^{(r)}$, $H_{1,1,1,2}^{(r)}$, $H_{1,1,1,3}^{(r)}$, $H_{1,1,1,4}^{(r)}$, which are the $r$-graphs extending from the 4-graphs shown in Fig. \ref{4-uniform} by $r-4$ times.	
\end{enumerate}
\end{theorem}
\begin{theorem}\citep{Lu2014Connected}\label{3-uniform-e}
	Let $r\geq 4$ . If the spectral radius of a connected $r$-uniform hypergraph $H$ is equal to $(r-1)!\sqrt[r]{4}$, then $H$ must be one of the following graphs:
	\begin{enumerate}
		\item[(1)] $C_n^{(r)}$, $\widetilde{D}_n^{(r)}$, $G_{1,2:n-8:1,2}^{(r)}$, $\widetilde{BD}_n^{(r)}$, $C_2^{(r)}$, $S_4^{(r)}$, which are $r$-graphs extending from those 3-graphs shown in Fig. \ref{equal} by $r-3$ times.
		\item[(2)] $E_{2,2,2}^{(r)}$, $E_{1,3,3}^{(r)}$, $E_{1,2,5}^{(r)}$, $F_{2,3,4}^{(r)}$, $F_{2,2,7}^{(r)}$, $F_{1,5,6}^{(r)}$, $F_{1,4,8}^{(r)}$, $F_{1,3,14}^{(r)}$, $G_{1,1:0:1,4}^{(r)}$ and $G_{1,1:6:1,3}^{(r)}$, which  are the $r$-graphs extending from corresponding 3-graphs by $r-3$ times.
		\item[(3)] $H_{1,1,2,2}^{(r)}$, which extends $r-4$ times from the hypergraph $H_{1,1,2,2}^{(4)}$ shown in Fig. \ref{4-uniform}.
	\end{enumerate}
\end{theorem}

Let $\mathcal{G}^{(r)}$ be the set of all $r$-graphs with spectral radius at most $(r-1)!\sqrt[r]{4}$, that is the collection of $r$-graphs mentioned in Theorems \ref{3-uniform-e} and \ref{r-uniform-s}. It is easy to see that $F^{(r)}_{1,3,14}$ has 19 edges, which is the most number of edges over $\mathcal{G}^{(r)}\backslash\{A^{(r)}_n, F^{(r)}_{1,1,n-3}, G^{(r)}_{1,1:n-6:1,1}, $ $F^{(r)}_{1,2,n-4},$ $ E^{(r)}_{1,1,n-2}, G^{(r)}_{1,2:n-7:1,1}, BD^{(r)}_{n}, C_n^{(r)}, \widetilde{D}^{(r)}_n, $ $G^{(r)}_{1,2:n-8:1,2}, \widetilde{BD}^{(r)}_n\}$. Then that is not difficult to find  the $r$-graphs with small spectral radius over $\H^{(r)}_{n}$, when  $n\geq 20$. Therefore, we only need to  consider those seven connected $r$-graphs shown in Fig. $\ref{small}$. 
\section{The effect on the spectral radius of $r$-graphs by perturbation ($r\geq3$)}
 In order to compare the spectral radius of those $r$-graphs listed in Fig. \ref{small}, we are going to first study three kinds of perturbations in this section.  Before coming to our results, two sequences of functions are needed.	

\begin{lemma}\label{fi}
Let

$$ f_1(x)=x, \;\;\;
 f_i(x)=\frac{x}{1-f_{i-1}(x)}, \;\;\; \mbox{for $i\geq 2$};$$
 $$	g_1(x)=x, \;\;\;g_2(x)=\frac{x}{(1-x)^2}, \;\;\;	g_i=\frac{x}{1-g_{i-1}(x)}, \;\;\;\mbox{for $i\geq 3$}.$$
 Denote $a_i$ and $b_i$ be the real number satisfying $f_i(a_i)=1$ and $g_i(b_i)=1$, respectively. For any positive integer $i\geq 1$, we have
\begin{enumerate}
		\item[(1)] $f_i(x)$ and $g_i(x)$ is increasing with respect to $x\in (0,a_{i-1})$ and $(0,b_{i-1})$, respectively, where $a_0=b_0=1$. Moreover, $a_{i}>a_{i+1}$ and $b_{i}>b_{i+1}$.	
		\item[(2)] $f_{i+1}(x)>f_{i}(x)$ with respect to $x\in(0,a_{i})$ and $g_{i+1}(x)>g_{i}(x)$ with respect to $x\in (0,b_{i})$. 
        \item[(3)] Let $c_2$ be the real number between $(0,1)$ satisfying $1-4c_2+2c_2^2=0$ {\it i.e.}, $c_2\approx 0.2929$. Then $f_{i+2}(x)<g_i(x)$ with respect to $x\in(0, c_i)$, where $c_i=\min\{a_{i+1},b_{i-1},c_2\}$ and $i\geq 2$.
\end{enumerate}

\end{lemma}
\p We will first prove Item (1) for $f_i(x)$ by induction on $i\in Z^+$. As we know, $f_i(0)=0$, for all $i\geq1$. It suffices to show that $f'_i(x)>0$ with respect to $x\in (0,a_{i-1})$.

 When $i=1$ and $2$, $f'_1(x)=1>0$ and $f'_2(x)=\frac{1}{(1-x)^2}>0$. Therefore, both $f_1(x)$ and $f_2(x)$ are increasing with respect to $x\in (0,1)$. Meanwhile, $a_1=1$, $a_2=\frac{1}{2}$ and $a_2<a_1$.

 Suppose that $f'_i(x)>0$ with respect to $x\in (0,a_{i-1}) $ for all $i\leq k$. Then $a_{k}<a_{k-1}$. In fact, $f_{k }(a_k)=\frac{a_k}{1-f_{k-1}(a_k)}=1$ and $f_{k-1}(a_k)=1-a_{k}<1$. Then $a_k<a_{k-1}$ follows from $f_{k-1}(x)$ is increasing when $x\in (0,a_{k-2})$. Then $f'_{k}(x)>0$ and $f_k(x)<1$ when $x\in (0,a_k)$.
 Therefore,
  $$f'_{k+1}(x)=\frac{1}{1-f_{k}(x)}+\frac{xf'_{k}(x)}{(1-f_{k}(x))}>0\;\;\;\; \mbox{when $x\in (0,a_k)$}.$$
 Similarly, we can prove
 $$g'_{i+1}(x)=\frac{1}{1-g_{i}(x)}+\frac{xg'_{i}(x)}{(1-g_{i}(x))}>0,\;\;\;\;\mbox{when $x\in (0,b_i)$.}$$  And $b_i>b_{i+1}$.

Secondly, we will show Item (2) is true for $f_i(x)$ by  induction on $i$.

If $i=1$, then $f_2(x)-f_1(x)=\frac{x}{1-x}-x=\frac{x^2}{1-x}>0$, when $x\in (0,1)$.

Suppose that $f_{i}(x)>f_{i-1}(x)$ with respect to $x\in(0,a_{i-1})$. By Item (1) and the inductive hypothesis, we have $(1-f_i(x))(1-f_{i-1}(x))>0$ when $x\in (0,a_i)$.
It follows that
$$f_{i+1}(x)-f_i(x)=\frac{x(f_i(x)-f_{i-1}(x))}{(1-f_i(x))(1-f_{i-1}(x))}>0,$$ when $x\in (0,a_i)$. 

Similarly,  $$g_{i+1}(x)-g_{i}(x)>0,\;\;\;\;\mbox{when $x\in (0,b_i)$}. $$

Finally, we will prove the Item (3) in the same way.

If $i=2$, then
\begin{equation*}
\begin{array}{lll}
f_4(x)-g_2(x)&=&\dfrac{x(1-2x)}{1-3x+x^2}-\dfrac{x}{(1-x)^2}\\
&=&-\dfrac{x^2(1-4x+2x^2)}{(1-x)^2(1-3x+x^2)}\\
&<&0.
\end{array}
\end{equation*}
Therefore, the result holds for $x\in (0,c_2)$, where $c_2$ is the root of $1-4x+2x^2=0$ between $(0,1)$ {\it i.e.}, $c_2\approx 0.2929$. Moreover, $a_3\approx 0.38$ and $b_1=1$. Hence, $f_4(x)<g_2(x)$ when $x\in (0,c_2)$.

Suppose that $f_{i+1}(x)<g_{i-1}(x)$ with respect to $x\in (0,c_{i-1})$, where $c_{i-1}=\min\{a_{i}, b_{i-2},$ $ c_2\}$. According to Item (1), $f_{i+1}(x)<1$  when $x\in(0,a_{i+1})$ and $g_{i-1}(x)<1$ when $x\in (0,b_{i-1})$. Hence,
\begin{equation*}
\begin{array}{lll}
f_{i+2}(x)-g_i(x)&=&\dfrac{x}{1-f_{i+1}}-\dfrac{x}{1-g_{i-1}(x)}\\
&=&\dfrac{x(f_{i+1}(x)-g_{i-1}(x))}{(1-f_{i+1}(x))(1-g_{i-1}(x))}\\
&<&0,
\end{array}
\end{equation*}
when $0<x<\min\{a_{i+1},b_{i-1},c_2\}$. The result holds.\q

{\bf Remark. } Note that $$f_6(x)-1=\frac{(1 - 2 x)(1 - 4 x + 2 x^2)}{(-1 + 5 x - 6 x^2 + x^3)}$$ and 
$$g_4(x)-1=\frac{(1 - x) (1 - 4 x + 2 x^2)}{-1 + 4 x - 3 x^2 + x^3}.$$ According to the monotonicity of $f_i(x)$ and $g_i(x)$, it gets  $b_{i-1}\leq a_{i+1}\leq c_2$ if $i\geq 5$ and $c_2<a_{i+1}<b_{i-1}$ if $2\leq i\leq 4$. Therefore, we can rewrite $c_i$ as follows.
\begin{equation*}
c_i=\left\{
\begin{array}{ll}
b_{i-1},& \;\;\;\mbox{if}\;\;\; i\geq 5;\\
c_2,&\;\;\;\mbox{if} \;\;\;2\leq i\leq 4.
\end{array}
\right.
\end{equation*}

Theorems \ref{hkl} and \ref{H+kl} is related to the effect on spectral radius of a $r$-graph by grafting an edge from a path to anther.

\begin{figure}[htbp]
	\begin{minipage}{0.5\textwidth}
	\centering
	\includegraphics[width=0.9\textwidth]{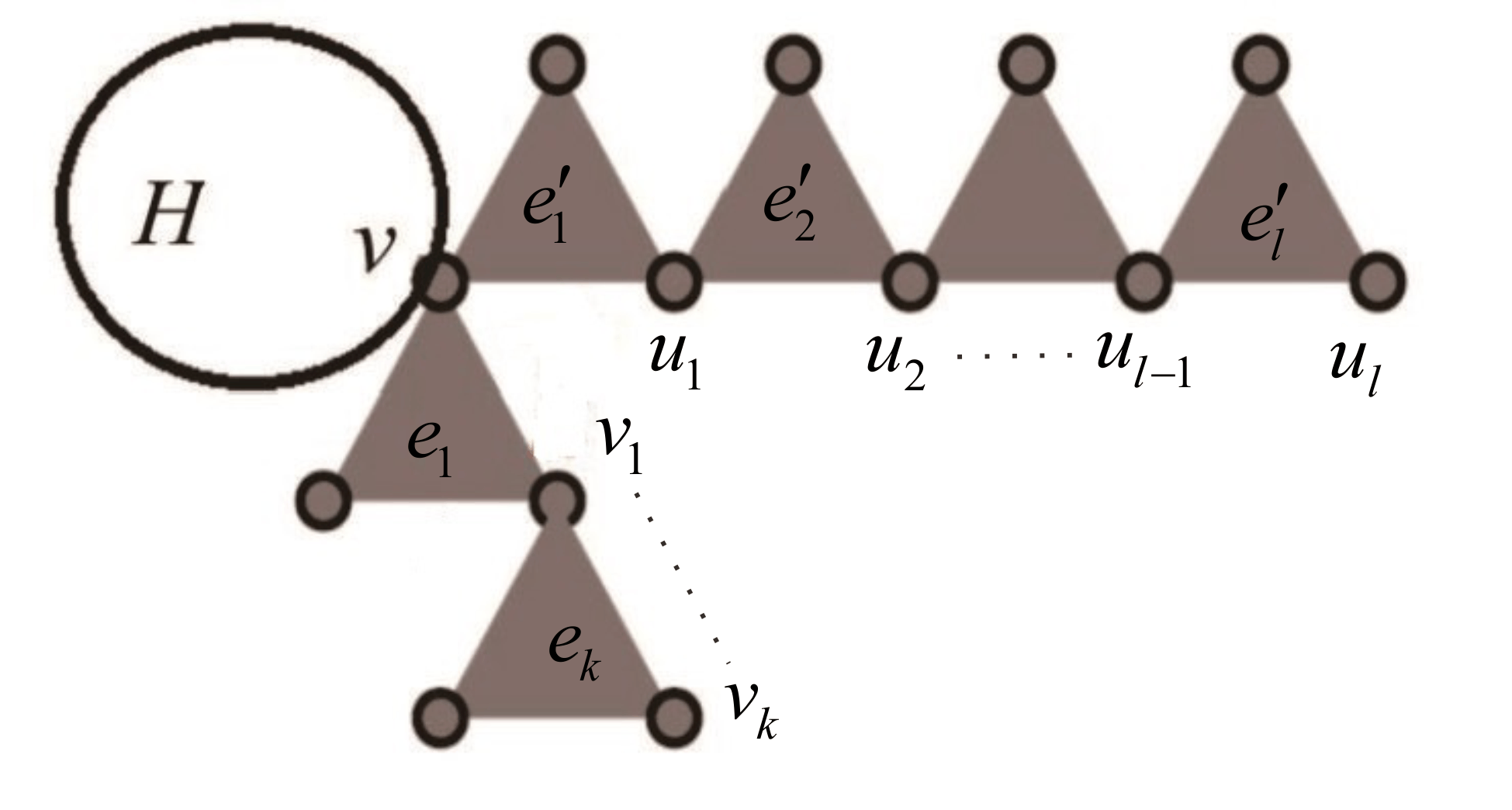}
	\caption{$H^{(r)}_{k,l}$}
	\label{hklp}
	\end{minipage}
	\begin{minipage}{0.5\textwidth}
		\centering
		\includegraphics[width=1.0\textwidth]{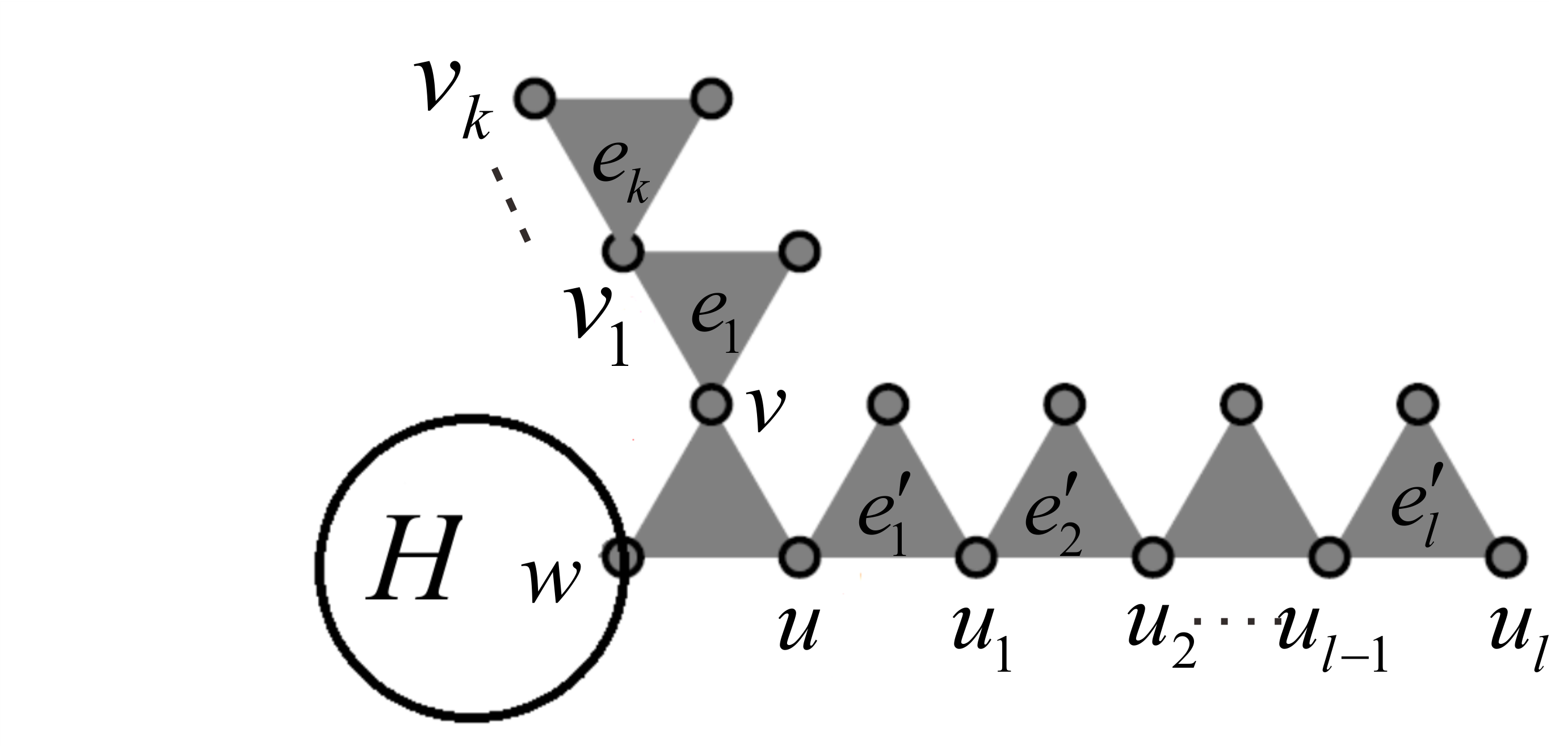}
		\caption{$\widetilde{H}^{(r)}_{k,l}$}
		\label{h+klp}
	\end{minipage}
\end{figure}

	Let $H$ be a nontrivial connected $r$-graph (not necessary to be simple) and $v$ be a vertex in $H$.  Suppose that $A_k^{(r)}=v_0e_1v_1e_2\cdots e_kv_k$ and $A_{l}^{(r)}=u_0e'_1u_1e'_2\cdots $ $e'_lu_l$, which are $r$-graphs extending from $A_k$ and $A_l$ by $r-2$ times. Denote $H^{(r)}_{k,l}$ as the $r$-graph obtained from $H$ by attaching $A_{k}^{(r)}$ and $A_{l}^{(r)}$ at $v$ such that $v_0=u_0=v$ (see Fig. \ref{hklp}).
	
	 Suppose $e=\{u,v,w\}$ is an edge in a connected $r$-graph $H$, where $u$ and $v$ are two leaf vertices and $w$ is a non-leaf vertex in $e$. Attaching $A_{k}^{(r)}$ and $A_{l}^{(r)}$ at $u$ and  $v$, respectively, such that $u=u_0$ and $v=v_0$, the resulting graph is $\widetilde{H}^{(r)}_{k,l}$ (see Fig. \ref{h+klp}).
\begin{theorem}\label{hkl}
	Let $H^{(r)}_{k,l}$ be the r-graph shown in Fig. \ref{hklp} and $l\geq k\geq 1$. Then $\rho(H^{(r)}_{k,l})>\rho(H^{(r)}_{k-1,l+1})$.
\end{theorem}			

\p Let $H^{(r)}_{k,l}$ be consistently $\alpha$-normal, where $\alpha=((r-1)!/\rho(H^{(r)}_{k,l}))^r$, and $B$ be the corresponding weighted incident matrix for $H^{(r)}_{k,l}$. In what follows, $f_i$ refers to $f_i(\alpha)$. As indicated in Fig.\ref{hklp}, $H^{(r)}_{k,l}$ has labeling with $B(u_{j-1},e'_{j})=f_{l+1-j}$ for $j=1,2,\cdots,l$, where $u_0=v_0=v$ and $B(v_{i-1},e_{i})=f_{k+1-i}$ for $i=1,2,\cdots,k$. Because $v$ is incident with at least 3 edges, we have $f_k+f_l<1$.

Define a weighted incident matrix $B'$ for $H^{(r)}_{k-1,l+1}$ with $B(u_{j-1},e'_{j})=f_{l+2-j}$ for $j=1,2,\cdots,l+1$ and $B(v_{i-1},e_{i})=f_{k-i}$ for $i=1,2,\cdots,k-1$. Otherwise, $B'(v,e)=B(v,e)$. Since $f_{k-1}=1-\frac{\alpha}{f_k}$ and $f_{l+1}=\frac{\alpha}{1-f_l}$, we have
\begin{align*}
&f_{k-1}+f_{l+1}-(f_k+f_l)\\
=&1-\frac{\alpha}{f_k}+\frac{\alpha}{1-f_l}-f_k-f_l\\
=&\frac{(f_k+f_l-1)(\alpha-f_k(1-f_l))}{(1-f_l)f_k}
\end{align*}
By Lemma \ref{fi} Item (2),		
$f_l\geq f_k>f_{k-1}=1-\frac{\alpha}{f_k}$. Then $\alpha-f_k(1-f_l)>0$.  Consequently, $f_{k-1}+f_{l+1}-(f_k+f_l)<0$ follows from $1-f_l>0$ and $f_k+f_l-1<0$, which implies $H_{k-1,l+1}$ is strictly $\alpha$-subnormal. Hence $\rho(H^{(r)}_{k,l})>\rho(H^{(r)}_{k-1,l+1})$ by Lemma \ref{sub-t}.\q
\begin{theorem}\label{H+kl}
 Let $\widetilde{H}^{(r)}_{k,l}$ be the r-graph shown in Fig.\ref{h+klp} and $x$ be the Perron-Frobenius vector, where $l\geq k\geq 1$. Then $x(u)\geq x(v)$ and $\rho(\widetilde{H}^{(r)}_{k,l})>\rho(\widetilde{H}^{(r)}_{k-1,l+1})$.
\end{theorem}		

\p Let $\widetilde{H}^{(r)}_{k,l}$ be consistently $\alpha$-normal, where $\alpha=((r-1)!/\rho(\widetilde{H}^{(r)}_{k,l}))^r$ and $B$ be the corresponding weighted incident matrix of $\widetilde{H}^{(r)}_{k,l}$. Then $\widetilde{H}^{(r)}_{k,l}$ has labeling with $B(u_{j-1},e'_{j})=f_{l+1-j}$ for $j=1,2,\cdots,l$ and $B(v_{i-1},e_{i})=f_{k+1-i}$, for $i=1,2,\cdots, k$, where $u_0=u$ and $v_0=v$. Furthermore, $B(v,e)=1-f_{k}$ and $B(u,e)=1-f_{l}$.

According to Item (2) of Lemma \ref{fi}, it gets $f_k\leq f_l$, {\it i.e.,} $B(v,e_1)\leq B(u,e'_1)$. Then $B(v,e)\geq B(u,e)$. Therefore, $x(v)\leq x(u)$ follows from Equation (\ref{eq1}).

Second, we consider the weighted incident matrix for $\widetilde{H}^{(r)}_{k-1,l+l}$, which has labeling with $B(v_{i-1},e_{i})=f_{k+1-i}$, for $i=1,2,\cdots, k-1$ and $B(u_{j-1},e'_{j})=f_{l+1-j}$ for $j=1,2,\cdots,l+1$, where $u_0=u$ and $v_0=v$.  $B(v,e)=1-f_{k-1}$ and $B(u,e)=1-f_{l+1}$. One can observe that
 $(1-f_k)(1-f_l)>\alpha$, due to $d(w)\geq 3$. Furthermore, $f_{l+1}>f_l\geq f_k$ follows from Lemma \ref{fi}. As a consequence,
\begin{align*}
&(1-f_{k-1})(1-f_{l+1})-(1-f_{k})(1-f_{l})	\\
=&(\frac{\alpha}{f_k})(1-\frac{\alpha}{1-f_l})-(1-f_{k})(1-f_{l})\\
=&\frac{1-f_l}{f_k}(\frac{\alpha}{1-f_l}-f_k)(1-f_k-\frac{\alpha}{1-f_l})\\
=&\frac{1-f_l}{f_k}(f_{l+1}-f_k)(1-f_k-\frac{\alpha}{1-f_l})\\
>&0, 
\end{align*}	
Then $\widetilde{H}^{(r)}_{k-1,l+l}$ is strictly $\alpha$-subnormal and the result follows by Lemma \ref{sub-t}.\q
\begin{corollary} Let $n$ be the number of edges in following $r$-graphs. Then 
	\begin{enumerate}
		\item[(1)] 	$\rho(F^{(r)}_{1,2,n-4})>\rho(F^{(r)}_{1,1,n-3})$ (for $n\geq 6$);
		\item[(2)]  $\rho(G^{(r)}_{1,1,n-6,1,1})>\rho(F^{(r)}_{1,1,n-3})$ (for $n\geq 6$);
		\item[(3)] $\rho(F^{(r)}_{1,4,n-1})>\rho(F^{(r)}_{1,3,n})$ (for $n\geq 5$);
		\item[(4)] $\rho(G^{(r)}_{1,2,n-7,1,1})>\rho(F^{(r)}_{1,2,n-4})>\rho(F_{1,1,n-3})$ (for $n\geq 7$).
	\end{enumerate}
  
\end{corollary}		
Combining Theorems \ref{hkl} and \ref{H+kl}, we can find the $r$-graph with smallest spectral radius among all connected $r$-graphs.	
\begin{corollary}
	$A^{(r)}_{n}$ is the  $r$-graph with minimal spectral radius over all connected $r$-graphs with size $n$.
\end{corollary}	

An edge $e$ is call a {\it 2-bridge} of a connected $r$-graph $H$ if $e$ contains exactly two non-leaf vertices and $H-e$ is disconnected. Let $u$ and $v$ be the two non-leaf vertex of the 2-bridge $e$. The contraction denoted by $H/e$ is obtained by deleting $e$ and identifying $u$ and $v$ to a new vertex $w$. The following theorem indicates that $\rho(H)$ will decrease after contracting a 2-bridge of $H$, if $\rho(H)< (r-1)!\sqrt[r]{4}$.

\begin{figure}[htbp]
	\centering
	\includegraphics[width=0.7\textwidth]{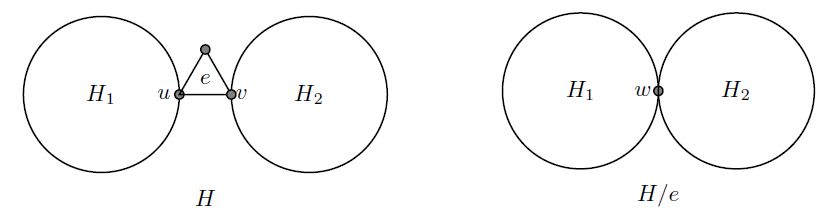}
	\caption{$H$ with a 2-bridge in Theorem \ref{divid-t}}
	\label{contrac}
\end{figure} 	
\begin{theorem}\label{divid-t}
Let $H$ be a connected r-graph and $e$ be a 2-bridge of $H$, where $u$ and $v$ are the two non-leaf vertices in $e$. Suppose that $H-e=H_1\cup H_2$ and $u\in V(H_1)$, $v\in V(H_2)$ (see Fig. \ref{contrac}). If $H$ is consistently $\alpha$-normal and $\alpha>\frac{1}{4}$, then $\rho(H/e)<\rho(H)$. 
\end{theorem}
\p Let $B$ be the corresponding weighted incident matrix of $H$, which is consistently $\alpha$-normal.  Suppose $x=B(v,e)$ and $y=B(u,e)$. Then $4\alpha=4xy\leq (x+y)^2$.  Define a weigthed incident  matrix $B'$ for $H/e$ with $B'(w,e)=B(u,e)$ if $e\in E(H_1)$ and $B(w,e)=B(v,e)$ if $e\in E(H_2)$. Otherwise, $B'(v,e)=B(v,e)$. Then 
$$\sum_{e\in E(H/e)}B(w,e)=2-(x+y)\leq 2-\sqrt{4\alpha}<1.$$
And $H/e$ is strictly  $\alpha$-subnormal. By Lemma \ref{sub-t}, $\rho(H/e)<\rho(H)$.
\q

According to Theorem \ref{divid-t}, one can get the following result.
\begin{corollary} Let $i,j,k,l,m,n$ be  nonnegative integers. Then
  $\rho(G_{i,j:k:l,m})<\rho(G_{i,j:k+1:l,m})$ and  $\rho(BD_n)<\rho(BD_{n+1})$. 
\end{corollary}
\section{The ordering of $r$-graphs by the spectral radius over $\H^{(r)}_n$}

In this section, we will compare the spectral radius of those seven hypergraphs listed in Fig.\ref{small} and then give the ordering of $r$-graphs with small spectral radius over $\H^{(r)}_n$, when $n\geq 20$.
\begin{theorem}
	$\rho(G^{(r)}_{1,1:n-6:1,1})\leq\rho(F^{(r)}_{1,2,n-4})$, when $n\geq 10$. The equality holds if and only if $n=10$.
\end{theorem}

	 \p  Let $G^{(r)}_{1,1:n-6:1,1}$ be consistently $\alpha_{n}$-normal and $B$ be the corresponding weighted matrix. For example, $G^{(r)}_{1,1:4:1,1}$ is consistently $\alpha_{10}$-normal.  As shown in Fig.\ref{small}, we have $B(v_{n-3},e_{n-2})=B(v_{n-2},e
	  _{n-1})=B(v_{n-1},e_n)=\alpha_{n}$,
	  $B(v_{n-3},e_{n-3})=B(v_{n-2},e
	  _{n-3})=B(v_{n-1},e_2)$ $=1-\alpha_{n}$, $B(v_{n-4},e_{n-3})=\frac{\alpha_{n}}{(1-\alpha_{n})^2}$. Meanwhile, $B(v_i,e_i)=g_i(\alpha_{n})$, $B(v_i,e_{i+1})=1-g_i(\alpha_{n})$, for $i=1,2,\cdots,n-4.$ Then $B(v_{n-4},e_{n-4})+B(v_{n-4},e_{n-3})=g_{n-4}(\alpha_{n})+\frac{\alpha_{n}}{(1-\alpha_{n})^2}=1$. 	 
	 
	Define a weighted incident matrix $B'$ for $F^{(r)}_{1,2,n-4}$ as follows.\\  $B'(v_{n-2},e_{n-1})=B'(v_{n-1},e_n)=\alpha_{n}$, $B'(v_{n-2},e_{n-2})=B'(v_{n-1},e_{n-3})=1-\alpha_{n}$, 
	$B'(v_{n-3},e_{n-2})=\frac{\alpha_{n}}{1-\alpha_{n}}$, $B'(v_{n-3},e_{n-3})=1-\frac{\alpha_{n}}{1-\alpha_{n}}=\frac{1-2\alpha_{n}}{1-\alpha_{n}}$, $B'(v_{n-4},e_{n-3})=\frac{\alpha_{n}}{1-2\alpha_{n}}$.
 Moreover,$B'(v_i,e_i)=f_{i}(\alpha_{n})$ and $B'(v_i,e_{i+1})=1-f_{i}(\alpha_{n})$, for $i=1,2,\cdots, n-4$. Hence, we need to check the value of  $B'(v_{n-4},e_{n-4})+B'(v_{n-4},e_{n-3})$ $=f_{n-4}(
 \alpha_{n})+\frac{\alpha_{n}}{1-2\alpha_{n}}$.
 
 First, we will show that $F^{(r)}_{1,2,6}$ is also consistently $\alpha_{10}$-normal. It is known that $\alpha_{10}$ satisfies the equation $g_6(x)-(1-\frac{x}{(1-x)^2})=0$, {\it i.e., }
$$\frac{(-1 + 4x - 3x^2 +x^3) (-1 + 6x - 9 x^2 + 3x^3)}{(-1 +x)^2 (1 - 6x + 10x^2 - 5x^3 +x^4)}=0.$$	
Moreover, $\alpha_{10}<\frac{1}{3}$ because $1-\alpha_{10}/(1-2\alpha_{10})>0$. Then $\alpha_{10}$ is a real number satisfying $-1 + 6\alpha_{10} - 9 \alpha_{10}^2 + 3\alpha_{10}^3=0$ and $\alpha_{10}\approx0.257773$.

On the other hand, 
simplifying $f_6(x)-(1-\frac{x}{(1-2x)})$, we have
$$f_6(x)-(1-\frac{x}{1-2x})=\frac{(-1 + 3x) (-1 + 6x - 9x^2 + 3x^3)}{(-1 + 2x) (-1 + 5x - 6x^2 + x^3)}.$$
It is easy to see that $f_6(\alpha_{10})+\frac{\alpha_{10}}{1-2\alpha_{10}}=1$. Therefore, $F^{(r)}_{1,2,6}$ are consistently $\alpha_{10}$-normal and  $\rho(G^{(r)}_{1,1:4:1,1})=\rho(F^{(r)}_{1,2,6})$.

Next, we will show that  $\rho(G^{(r)}_{1,1:n-6:1,1})<\rho(F^{(r)}_{1,2,n-4})$ for $n\geq 11$.  
It is equivalent to prove that $f_{n-4}(\alpha_{n})+\frac{\alpha_{n}}{1-2\alpha_{n}}>1$, when $g_{n-4}(\alpha_{n})+\frac{\alpha_{n}}{(1-\alpha_{n})^2}=1$. 

Let $$L_i(x)=f_i(x)+\frac{x}{1-2x}-(g_i(x)+\frac{x}{(1-x)^2}), $$ where $i\geq 7$. We are going to show that $L_i(x)\geq 0$ with respect to $x\in (0.25,\alpha_{i+4})$, proceed by induction on $i\geq 7$.

When $i=7$ or $8$,
$$f_{7}(x)=\frac{x (-1 + 5x - 6x^2 + x^3)}{-1 + 6x - 10x^2 + 4x^3};$$
$$g_{7}(x)=\frac{x (1 - 6 x + 10x^2 - 5x^3 + x^4)}{1 - 7x + 15 x^2 - 11 x^3 + 3 x^4};$$
$$f_{8}(x)=\frac{x - 6 x^2 + 10 x^3 - 4 x^4}{1 - 7x + 15 x^2 - 10 x^3 + x^4};$$
$$g_{8}(x)=\frac{-x (1 - 7x + 15 x^2 - 11 x^3 + 3 x^4)}{-1 + 8 x - 21x^2 + 21 x^3 - 8 x^4 + x^5}.$$
Solving the equation $g_{7}(x)+\frac{x}{(1-x)^2}-1=0$, we can get $\alpha_{11}\approx 0.25672<0.2568$. Moreover,
$$
L_7(x)=-\frac{x^3 (1 - 11 x + 45 x^2 - 85 x^3 + 76 x^4 - 31 x^5 + 3 x^6 + 
	x^7)}{(-1 + x)^2 (-1 + 2 x) (1 - 4 x + 2 x^2) (1 - 7 x + 
	15 x^2 - 11 x^3 + 3 x^4)}>0,$$
when $x\in (0.25,0.2568)$.  Then $L_7(x)>0$ with respect to $x\in(0.25,\alpha_{11})$ and $F_{1,2,5}$ is strictly and consistently $\alpha_{11}$-supernormal.

Similarly, we have $\alpha_{12}\approx 0.255903<0.256$. Moreover, when $x\in (0.25,0.256)$, 
$$L_8(x)=\frac{x^3 (1 - 11 x + 46 x^2 - 93 x^3 + 97 x^4 - 54 x^5 + 18 x^6 - 3 x^7)}{(-1 + x)^2 (-1 + 2 x) (-1 + 5 x - 6 x^2 + x^3) (1 - 6 x + 10 x^2 - 
	5 x^3 + x^4)}>0.$$
 Then $L_8(x)>0$ with respect to $x\in(0.25,\alpha_{12})$. The result holds.

Suppose $L_{i-1}(x)> 0$ with respect to $x\in (0.25,\alpha_{i+3})$. In what follows, we will show the result is true for $i$.

Since $g_{i}(\alpha_{i+4})+\frac{\alpha_{i+4}}{(1-\alpha_{i+4})^2}=1$ and $g_{i}(x)+\frac{x}{(1-x)^2}$ is increasing with respect to $x$, it gets $$g_{i}(x)+\frac{x}{(1-x)^2}=\frac{x}{1-g_{i-1}(x)}+\frac{x}{(1-x)^2}<1, $$ when $x\in (0.25,\alpha_{i+4})$, {\it i.e.,}
$x<(1-\frac{x}{(1-x)^2})(1-g_{i-1}(x)).$
Meanwhile, 
$$f_8(x)-\frac{x}{(1-x)^2}=\frac{-x^2 (-1 + 6 x - 9 x^2 + 3 x^3)}{(-1 + x)^2 (-1 + 5 x - 6 x^2 + x^3)}>0, $$for any $x\in (0.25,0.2568)$. Therefore, $f_i(x)>\frac{x}{(1-x)^2}$ for any integer $i\geq 8$.
Hence, when $i\geq 9$ and $x\in (0.25,\alpha_{i+4})\subset (0.25,\alpha_{i+3})$, we have
\begin{equation*}
\begin{array}{lll}
L_i(x)&=&f_i(x)+\dfrac{x}{1-2x}-(g_i(x)+\dfrac{x}{(1-x)^2})\\
&=&\dfrac{x}{1-f_{i-1}(x)}+\dfrac{x}{1-2x}-(\dfrac{x}{1-g_{i-1}(x)}+\dfrac{x}{(1-x)^2})\\
&=&\dfrac{x(f_{i-1}(x)-g_{i-1}(x))}{(1-f_{i-1}(x))(1-g_{i-1}(x))}+\frac{x}{1-2x}-\frac{x}{(1-x)^2}\\
&\geq& \dfrac{(1-\frac{x}{(1-x)^2})(f_{i-1}(x)-g_{i-1}(x))}{(1-f_{i-1}(x))}+\frac{x}{1-2x}-\frac{x}{(1-x)^2} \\
&\geq &f_{i-1}(x)-g_{i-1}(x)+\frac{x}{1-2x}-\frac{x}{(1-x)^2}\\
&>& 0.
\end{array}
\end{equation*}
In summary, $\rho(G^{(r)}_{1,1:n-6:1,1})\leq \rho(F^{(r)}_{1,2,n-4})$, when $n\geq10$.\q

{\bf Remark.} Actually, $\rho(G^{(r)}_{1,1:n-6:1,1})>\rho(F^{(r)}_{1,2,n-4})$, for $6\leq n\leq 9.$ By calculating, their spectral radius, when $r=3$, are listed in the following table and the cases for $r\geq 4$ follows from Lemmas \ref{same} and \ref{lem1}.
\begin{center}
	\begin{tabular}{c|c|c}\hline
		$n$& $\rho(G^{(3)}_{1,1:n-6:1,1})$ & $\rho(F^{(3)}_{1,2,n-4})$ \\ \hline\hline
		6 & 3.1023  &3.0703 \\\hline
		7& 3.1188 & 3.1023\\\hline
	8& 3.1295&3.1215\\\hline
		9& 3.1370 & 3.1340\\\hline
	\end{tabular}
\end{center}
\begin{theorem}\label{th2}
	$\rho(G^{(r)}_{1,2:n-7:1,1})<\rho(BD^{(r)}_n)$, $\rho(F^{(r)}_{1,2,n-4})<\rho(E^{(r)}_{1,1,n-2})$, for $n\geq 7$.	
\end{theorem}
\p We firstly prove $\rho(F^{(r)}_{1,2,n-4})<\rho(E^{(r)}_{1,1,n-2})$ for $n\geq 7$.
Suppose that  $E^{(r)}_{1,1,n-2}$ is consistently $\alpha$-normal and $B$ is the corresponding weighted incident matrix. As shown in Fig. \ref{small}, it is labeled with $B(v_{n-2},e_{n-1})=B(v_{n-2},e_{n})=\alpha$, $B(v_{n-2},e_{n-2})=1-2\alpha$, $B(v_{n-3},e_{n-2})=\frac{\alpha}{1-2\alpha}$, $B(v_{n-3},e_{n-3})=1-\frac{\alpha}{1-2\alpha}=\frac{1-3\alpha}{1-2\alpha}$ and $B(v_{n-4},e_{n-3})=\frac{\alpha(1-2\alpha)}{1-3\alpha}$. Meanwhile, $B(v_i,e_i)=f_{i}(\alpha)$, for $i=1,2,\cdots,n-4$.
 Therefore, $B(v_{n-4},e_{n-4})+B(v_{n-4},e_{n-3})=f_{n-4}(\alpha)+\frac{\alpha(1-2\alpha)}{1-3\alpha}=1$.

Define a weighted incident matrix $B'$ for $F^{(r)}_{1,2,n-4}$ as follows. $B'(v_{n-2},e_{n-1})=B'(v_{n-1},e_{n})$ $=\alpha$, $B'(v_{n-2},e_{n-2})=B'(v_{n-1},e_{n-3})=1-\alpha$, $B'(v_{n-3},e_{n-2})=\frac{\alpha}{1-\alpha}$, $B'(v_{n-3},e_{n-3})=1-\frac{\alpha}{1-\alpha}=\frac{1-2\alpha}{1-\alpha}$. And $B'(v_{n-4},e_{n-3})=\frac{\alpha}{1-2\alpha}$. Otherwise, $B'(v,e)=B(v,e)$.

Since $\alpha>\frac{1}{4}$, we have
$$\frac{\alpha(1-2\alpha)}{1-3\alpha}-\frac{\alpha}{1-2\alpha}>0.$$
Consequently,  $B'(v_{n-4},e_{n-3})+B'(v_{n-4},e_{n-4})<1$ and  $F^{(r)}_{1,2,n-4}$ is strictly $\alpha$-subnormal. By Theorem \ref{sub-t}, it gives $\rho(F^{(r)}_{1,2,n-4})<\rho(E^{(r)}_{1,1,n-2})$, for $n\geq 7$.

That means $\rho(F^{(r)}_{1,2,n-4})$ will increase by grafting $e_n$ from $v_{n-1}$ to $v_{n-2}$. Therefore, $\rho(G^{(r)}_{1,2:n-7:1,1})<\rho(BD^{(r)}_n)$ will directly follow in the same way.\q

\begin{theorem}
	$\rho(G^{(r)}_{1,2:n-7:1,1})>\rho(E^{(r)}_{1,1,n-2})$, for $n\geq 7$.
\end{theorem}	
	
\p When $n=7$, by calculating we have $\rho(G^{(3)}_{1,2:0:1,1})\approx3.1426$ and $\rho(E^{(3)}_{1,1,5})\approx 3.1215$. According to Lemmas \ref{same} and \ref{lem1}, we have $\rho(G^{(r)}_{1,2:0:1,1})>\rho(E^{(r)}_{1,1,5})$, for any integer $r\geq 3$. 

Next, we will show the result holds for $n\geq 8$. Let $G^{(r)}_{1,2:n-7:1,1}$ be consistently $\alpha$-normal and $B$ be the corresponding weighted incident matrix. It is labeled with $B(v_i,e_i)=g_i(\alpha)$, for $i=1,2,\cdots,n-5$. And $B(v_{n-3},e_{n-2})=B(v_{n-1},e_{n})=B(v_{n-2},e_{n-1})=\alpha$. Then $B(v_{n-3},e_{n-3})=B(v_{n-2},e_{n-4})=B(v_{n-1},e_2)=1-\alpha$, $B(v_{n-4},e_{n-3})=\frac{\alpha}{1-\alpha}$ and $B(v_{n-4},e_{n-4})=1-\frac{\alpha}{1-\alpha}=\frac{1-2\alpha}{1-\alpha}$, $B(v_{n-5},e_{n-4})=\frac{\alpha}{1-2\alpha}$. Since $G_{1,2:n-7:1,1}$ is consistently $\alpha$-normal, $g_{n-5}(\alpha)+\frac{\alpha}{1-2\alpha}=1$, {\it i.e.}, $\frac{\alpha}{1-g_{n-6}(\alpha)}=1-\frac{\alpha}{1-2\alpha}=\frac{1-3\alpha}{1-2\alpha}$. 

Define a weighted incident matrix $B'$ for $E^{(r)}_{1,1,n-2}$ as follows.\\ $B'(v_{n-2},e_{n-1})= B'(v_{n-2},e_{n-1}) =\alpha$, $B'(v_{n-2},e_{n-2})=1-2\alpha$, $B'(v_{n-3},e_{n-2})=\frac{\alpha}{1-2\alpha}$, $B'(v_{n-3},e_{n-3})=1-\frac{\alpha}{1-2\alpha}=\frac{1-3\alpha}{1-2\alpha}$ and $B'(v_{n-4},e_{n-3})=\frac{\alpha(1-2\alpha)}{1-3\alpha}$. Moreover, $B'(v_i,e_{i})=f_{i}(\alpha)$ and $B'(v_i,e_{i+1})=1-f_{i}(\alpha)$, for $i=1,2,\cdots,n-4$. 

Basically, $\alpha\in(0,c_{n-6})$. In fact, by Theorem \ref{divid-t} $\rho(G^{(r)}_{1,2:n-7:1,1})>\rho(G^{(r)}_{1,2:0:1,1})>3.1426$ and $\alpha<0.2578<c_2$. Moreover, $g_{n-7}(\alpha)<g_{n-5}(\alpha)<1$, then $\alpha<b_{n-7}$. Furthermore,
by Lemma \ref{fi} Item (3), we have
\begin{equation*}
\begin{array}{lll}
f_{n-4}(\alpha)+\dfrac{\alpha(1-2\alpha)}{1-3\alpha}-1&=&f_{n-4}(\alpha)+(1-g_{n-6}(\alpha))-1\\
&=&f_{n-4}(\alpha)-g_{n-6}(\alpha)\\
&<&0.
\end{array}
\end{equation*}	
Consequently, $B'(v_{n-4},e_{n-4})+B'(v_{n-4},e_{n-3})<1$. Then $E^{(r)}_{1,1,n-2}$ is $\alpha$-subnormal and $\rho(G^{(r)}_{1,2:n-7:1,1})>\rho(E^{(r)}_{1,1,n-2})$.\q

	 In conclusion, we derive the following.
	\begin{theorem}
 Let  $G$ be an r-graph in $\mathcal{H}^{(r)}_n\backslash\{A^{(r)}_n, F^{(r)}_{1,1,n-3}, G^{(r)}_{1,1:n-6:1,1}, $ $F^{(r)}_{1,2,n-4},$ $ E^{(r)}_{1,1,n-2}, G^{(r)}_{1,2:n-7:1,1}, BD^{(r)}_{n}, C_n^{(r)}, \widetilde{D}^{(r)}_n, G^{(r)}_{1,2:n-8:1,2}, \widetilde{BD}^{(r)}_n\}$ and $n\geq 20$. Then $\rho(G)>(r-1)!\sqrt[r]{4}$. Moreover,  
		$\rho(A^{(r)}_n)<\rho(F^{(r)}_{1,1,n-3})<\rho(G^{(r)}_{1,1:n-6:1,1})<\rho(F^{(r)}_{1,2,n-4})<\rho(E^{(r)}_{1,1,n-2})<\rho(G^{(r)}_{1,2:n-7:1,1})<\rho(BD^{(r)}_{n})<\rho(C_n^{(r)})=\rho(\widetilde{D}^{(r)}_n)=\rho(G^{(r)}_{1,2:n-8:1,2})=\rho(\widetilde{BD}^{(r)}_n)=(r-1)!\sqrt[r]{4}$.	
	\end{theorem}
\vspace{1cm}
{\bf References}	
\bibliographystyle{plain}
\bibliography{reference2}
\end{document}